\documentclass[12pt]{article}
\usepackage{amsmath,amssymb}

\newtheorem{thm}{Theorem}
\newtheorem{lemma}{Lemma}
\newtheorem{prop}{Proposition} 
\newtheorem{rem}{Remark}

\newcommand{\mb}{\mbox}
\newcommand{\x}{\times}
\newcommand{\g}{\mathfrak g}

\newcommand{\ot}{\otimes}

\newcommand{\<}{\langle}
\renewcommand{\>}{\rangle}
\renewcommand{\a}{\alpha}
\renewcommand{\b}{\beta}

\newcommand{\C}{\mathbb C}
\newcommand{\ds}{\displaystyle}

\newcommand{\D}{\Delta}
\newcommand{\e}{\varepsilon}
\newcommand{\alg}{algebra}
\newcommand{\ij}{a_{ij}}
\newcommand{\ai}{\alpha_i}
\newcommand{\Spec}{\mbox{Spec }}
\newcommand{\R}{\mathcal R}
\newcommand{\mor}{morphism}
\newcommand{\U}{\mathcal U}

\begin{document}
\setlength{\baselineskip}{16pt}

\title{Invariants of links with flat connections in their complements.II.
Holonomy $R$-matrices related to quantized universal enveloping
algebras at roots of 1.}

\author{R. Kashaev,  N.Reshetikhin}
\maketitle

\begin{abstract}
Holonomy $R$-matrices parametrized by finite-dimensional representations
are constructed for quantized
universal enveloping algebras of simple Lie algebras at roots of 1.
\end{abstract}

\section{Introduction}

In the previous paper \cite{KR} we gave the construction of
invariants of tangles with flat $G$-connections in the complement.
The construction was based on the notion of holonomy $R$-matrices.
Such $R$-matrices are operator-valued functions on $G\times G$
which satisfy the holonomy Yang-Baxter equation \cite{KR}.

In this paper we show that quantized universal enveloping algebras
at roots of 1 provide examples of such $R$-matrices.

In section 2 we present basic facts about quantized universal
enveloping algebra when the quantization parameter is generic. In
section 3 we show that the regular part of the universal $R$-
matrix specialized at roots of 1 give examples of holonomy $R$-
matrices. In section 4 we analyze the invariants of tangles with a
flat connection in the complement derived from such holonomy
$R$-matrices. This work was supported by the NSF grant DMS-0070931
and by the DMS grant RM1-2244.

\section{Quantized Universal Enveloping Algebras}

Let ${\g}$ be a simple Lie \alg \ of rank $r$ with the root system
$\D$. Denote by $P$ its weight lattice and by $Q$ its root
lattice. Fix simple roots $\a_1,\dots ,\a_r\in\D$ and denote by
$(\ij)^r_{i,j=1}$ the corresponding Cartan matrix. Denote by
$\omega_1,\dots ,\omega_r\in P$ the basis of fundamental weights (dual to
the basis of simple roots
$\{\ai\}$) and let $d_i$ be the the length of the $i$-th simple
root.

\subsection{Quantized universal enveloping algebras}
Let $Q\leq M\leq P$ be a lattice. The quantized
universal enveloping \alg \ $U^M_q({\g})$ is the associative \alg
\ with 1 over ${\C}(q)$ generated by $L_\mu$, $\mu\in M$, and
$E_i,F_i$, $i=1,\dots ,r$ with defining relations:
\begin{eqnarray*}
L_\mu L_\nu &=& L_\nu L_\mu \ , \quad
L_0=1 \ , \\ L_\mu E_i &= & q^{\a_i(\mu)}E_iL_\mu \\
L_\mu F_i &=& q^{-\ai (\mu)}F_iL_\mu  \ , \qquad
E_iF_j -F_jE_i = \delta_{ij}
\frac{L_{\ai}-L_{\ai}^{-1}}{q_i-q_i^{-1}} \ , \\
&&\sum^{1-\ij}_{k=0} (-1)^k
\left[ \begin{array}{c} 1-\ij\\ k\end{array}\right]_{q_i}\,
E_i^{1-k-\ij} E_j\, E_k^k =0 \ , \quad i\neq j\\
&&\sum^{1-\ij}_{k=0} (-1)^k
\left[ \begin{array}{c} 1-\ij\\ k\end{array}\right]_{q_i}\,
F_i^{1-k-\ij} F_j\, F_k^k =0 \ , \quad i\neq j
\end{eqnarray*}
Here \ $q_i=q^{d_i}$,
\[
\left[ \begin{array}{c} m\\n\end{array}\right]_q\,
= \frac{[m]_q!}{[m-n]_q![n]_q!} \ ,\qquad
[n]_q!=[n]_q\dots [2]_q[1]_q \ , \qquad
\left[n\right]_q = \frac{q^n-q^{-n}}{q-q^{-1}} \ .
\]
The map $\D$ acting on generators as
\begin{eqnarray}\label{comult}
\D L_\mu &=& L_\mu\ot L_\mu \ , \\
\D E_i &=& E_i\ot 1 + L_{\ai}\ot E_i \ ,\\
\D F_i &=& F_i\ot L_{-\ai} + 1\ot F_i
\end{eqnarray}
extends to the homomorphism of \alg s $\D:U^M_q({\g})\to
U^M_q({\g})\ot U^M_q({\g})$. The pair $(U^M_q({\g}),\D)$ is
a Hopf \alg \ with the counit $\e(L_\mu)=1$, $\e(E_i)=\e(F_i)=0$.

It is clear that if $P\leq M\leq M'\leq Q$ we have the embedding
$U^M_q({\g})\hookrightarrow  U^{M'}_q({\g})$ of Hopf \alg s.

\subsection{Quantum Weyl group}
Let $B_W$ be the braid group associated with the Weyl group $W$,
\[
B_W=\{\mb{generated by} \ T_i\mid \mb{with defining relations} \
\underbrace{T_iT_jT_iT_j\dots}_{m_{ij}} =
\underbrace{T_jT_iT_jT_i\dots}_{m_{ij}}\}
\]
Here $m_{ij}=2$ if $i$ and $j$ are not connected in the Dynkin
diagram, $m_{ij}=3$ if $a_{ij}a_{ji}=1$, $m_{ij}=4$ if $a_{ij}a_{ji}=2$ and
$m_{ij}=6$ if $a_{ij}a_{ji}=3$.

This group acts on $U^M_q({\g})$ by auto\mor s $[L]$:
\begin{eqnarray*}
T_i(L_\mu) &=& L_{s_i(\mu)} \ , \\
T_i(E_i) &=& -F_iL_{\ai} \ ,\\
T_i(E_j) &=&\sum_{r+s=-\ij}(-1)^r q^s_i E^{(r)}_i E_j E_i^{(s)}\\
T_i(F_i) &=& -L_{\ai}^{-1} E_i \\
T_i(F_j) &=& \sum_{r+s=-\ij} (-1)^r q^s_i F_i^{(s)} F_jF_i^{(r)}
\end{eqnarray*}
where \ $X_i^{(r)}=\frac{X_i^r}{[r]_{q_i}!}$.

Fix a reduced decomposition
of the longest element $w_0\in W$. If $w_0=s_{j_1}\dots s_{j_N}$,
$N=|\D_+|$, then
\begin{equation}\label{order}
\b_a = s_{j_1}\dots s_{j_{a-1}}\a_{j_a} \ , \ a=1,\dots, N \ ,
\end{equation}
were $\a_1\dots \a_r$ are simple roots. This gives a total convex
ordering $\beta_1<\dots<\beta_N$ on the set of roots $\Delta$ of $\mathfrak g$. Such
construction give all convex orderings and vice versa.

According to \cite{Lu} define root elements of $U_q({\mathfrak g})$ as
\[
E_{\beta_a}=T_{j_1}\dots T_{j_{a-1}} (E_{j_a}),
\]
\[
F_{\beta_a}=T_{j_1}\dots T_{j_{a-1}} (F_{j_a})
\]

The elements
\[
{E}^{k_1}_{\b_1}\dots { E}^{k_N}_{\b_N} \ L_\mu \, {
F}^{\ell_N}_{\b_N}\dots {F}_{\b_1}^{\ell_1}
\]
form a linear basis in the algebra $U_q({\mathfrak g})$.

\subsection{Integral form of quantized universal eneveloping algebras}
Define the ${\C}[q,q^{-1}]$-sub\alg \ ${\mathcal
U}^M_q({\g})\subset U^M_q({\g})$ as the smallest $B_W$-stable
${\C}[q,q^{-1}]$-sub\alg \ of $U^M_q({\g})$ containing the
elements
\[
\overline{E_i} = (q_i-q_i^{-1})E_i \ , \qquad
\overline{F_i} = (q_i-q_i^{-1})F_i \ .
\]
Set \ $\overline{E_\a}=(q_\a-q_\a^{-1})E_\a$, \
$\overline{F_\a}=(q_\a-q_\a^{-1})F_\a$, then monomials
\begin{equation}\label{pbw}
{\bar E}^{k_1}_{\b_1}\dots {\bar E}^{k_N}_{\b_N} \
L_\mu \, {\bar F}^{\ell_N}_{\b_N}\dots {\bar F}_{\b_1}^{\ell_1}
\end{equation}
form a linear basis in $\U^M_q({\g})$. Here we used the
enumeration of positive roots corresponding to a
reduced decomposition of the longest element of the Weyl group
(see above).

\subsection{Poisson Lie groups $G$ and $G^*$} \label{PL}

It is well known that the algebra ${\mathcal
U}^Q_q({\g})$ can be regarded as a Hopf algebra deformation of the
algebra of polynomial functions on the Poisson Lie group
$G^*=\{(b_+,b_-)\in B_+\times B_-| [b_+]_0=[b_-]_0^{-1} \}$.
Notice that as a Lie group $G^*$ is naturally isomorpic to the
semidirect product $H\ltimes (N_+\x N_-)$ where $H$ act naturally
on $N_\pm$. The tangent Lie bialgebra for
this Poisson Lie group is dual to the standard Lie bialgebra
structure on $\mathfrak g$ \cite{CP}. In this sense the Poisson
Lie group $G^*$ is dual to the Poisson Lie group $G$.

Similarly for any lattice $M, \ Q\leq M\leq P$ the covering group
$G_M^*$ of $G^*$ is also a Poisson Lie group which is dual to
the standard Poisson Lie group structure on $G$ in a sence that their
tagent Lie bialgebras are dual.

The algebra $C(G_M^*)$ of algebraic functions on the
Poisson Lie group $G_M^*$ is a
Poisson Hopf algebra. As a Poisson algebra it is generated by
elements $k_\mu, e_i, f_i$, $\mu\in M$, $i=1,\dots, r$
with defining relations
\[
\{k_\mu,k_\nu\}=0, \{k_\mu, e_j\}=\mu(\alpha_i)k_\mu e_j
\]
\[
\{e_i,f_j\}=\delta_{ij} (k_{\alpha_i}-k_{\alpha_i}^{-1})
\]
\[
 \underbrace{\{e_i,\dots,\{e_i }_{-a_{ij}+1} ,
 e_j\}^{(d_ia_{ij})}\dots \}^{(-d_ia_{ij})} = 0
\]
\[
 \underbrace{\{f_i,\dots,\{f_i }_{-a_{ij}+1} ,
 f_j\}^{(d_ia_{ij})}\dots \}^{(-d_ia_{ij})} = 0
\]
where $\{ X,Y \}^{(n)} = \{X,Y\}-nXY$.

The comultiplication acts on generators as in (\ref{comult}). The
elements $k_\mu$ are coordinate functions on the Cartan subgroup of
$G_M^*$ which is a finte cover of the Catran subgroup of $G$.
The elements $e_i$ and $f_i$
are coordinate functions on the nilpotent subgroups $N^\pm\subset
G^*$ corresponding to the simple roots.

The braid group $B_W$ acts on $C(G_M^*)$ by Poisson automorphisms.

\begin{eqnarray*}
\tau_i(k_\mu) &=& k_{s_i(\mu)} , \\
 \tau_i(e_i) &=&-f_i k_{\alpha_i}^{-1}, \\
 \tau_i(f_i) &=& -e_i k_{\alpha_i} , \\
 \tau_i(e_j) &=& \frac{(-1)^{a_{ij}}}{(-a_{ij})!}
  \{e_i,\dots\{e_i,e_j\}^{(a_{ij}d_i)}\}^{(d_i(a_{ij}+2))}\dots\}^
  {d_i(-a_{ij}-2)}, \\
 \tau_i(f_j) &=& \frac{1}{(-a_{ij})!} \{f_i,\dots \{f_i,f_j\}^{(a_{ij}d_i)}\}^
  {(d_i(a_{ij}+2))}\dots \}^{(d_i(-a_{ij}-2))}
\end{eqnarray*}

One can define coordinates corresponding to all positive and
negative roots on $G_M^*$ similarly to how it was done for
$U^M_q({\mathfrak g})$.

Fix a linear isomorphism between $\U_q^M({\mathfrak g})$
and $C(G_M^*)$ by identifying monomials (\ref{pbw}) with
corresponding monomials
in $k_\mu$, $e_\a$, $f_\a$. Then it is clear that the Hopf algebra structure
on $\U_q^M({\mathfrak g})$ is a Hopf algebra deformation
of the Poisson Hopf algbera structure on $C(G_M^*)$ described above.

\subsection{Symplectic leaves of $G_M^*$}
According to the general structural facts about Poison Lie groups
symplectic leaves of  $G^*$ are
orbits of the (local) dressing action of the dual Poisson
Lie group $G$ \cite{STS}. This action can be describe as follows.

Let $I: G^*\to G$ be the natural map $(x_+,x_-)\mapsto
x_+(x_-)^{-1}$. This map intertwines the dressing action with the
adjoint action of $G$ on $G$.

The map $I$ brings the dressing action of $G$ on $G^*$ to
the action of $G$ on itself by conjugations, i.e.
if we will write $g: x\mapsto g(x)$ for the dressing action of
$g\in G$ on $x\in G^*$ we have:
\[
I(g(x))=gI(x)g^{-1}
\]
Thus, orbits of dressing action in $G^*$ are connected components of
orbits of adjoint action of $G$ on itself.

The image of the map $I$ is open
dense in $G$. Over generic point in $G$ it is a branched
cover map with $2^r$-fibers and it
gives an isomorphism between a  neighborhood of 1 in $G$ and
neighborhoods of points $(\sigma, \sigma^{-1}) \in G^*$ where
$\sigma\in H$, $H$ is a Cartan subgroup in $G$ and $\sigma^2=1$. Using this ismorphisms we
can identify these neighborhoods of $G^*$ and  $G$.
We will call it a realization of  $G^*$ on
$G$.

The natural projection $G_M^*\to G^*$ is Poisson and
is a finite cover. Therefore
symplectic leaves of $G_M^*$ are connected components
of preimages of symplectic leaves in $G^*$.

\subsection{Formal Poisson Lie group $G_M^*$}\label{epsilon}
Let $\Gamma_M$ be the finite subgroup in $G_M^*$ which is the pre-image
of $1\in G^*$ with respect to the natural projection $G_M^*\to G^*$.
Let $\varepsilon_\mu\in \Gamma_M$ be the element corresponding to
the weight $\mu\in M/Q$.

Denote by $F(G_M^*)$ the completion of $C(G_M^*)$ by formal power
series in $k_\mu\epsilon_\mu^{-1}-1$, $e_\alpha$ and $f_\alpha$.
This Poisson Hopf
algebra is the formal Poisson Lie group $G_M^*$. Instead of formal
variables $k_\mu$ we can work with $z_\mu$ such that
$k_\mu=\varepsilon_\mu\exp(z_\mu)$.

\section{Quantized universal enveloping algebras at roots of 1}
\subsection{ The algebra $\U^M_\e({\g})$ and its center} Let $\ell$ be an odd integer such that $\ell >
\max_i(d_i)$ and $\e\in{\C}$ be a primitive $\ell$-th root of~1.
Denote by $\U^M_\e({\g})$ be the quotient \alg
\[
\U^M_\e({\g}) = \frac{\U^M_q({\g})}{(q-\e)\U^M_q({\g})} \ .
\]
The center $Z^M_\e =Z(\U^M_\e({\g}))$ has natural structure of Poisson
\alg \ and, as a Poisson \alg, it acts on $\U^M_\e({\g})$ by
derivations \cite{DKP}.

Denote by $Z^M_0$ the sub\alg \ in $Z^M_\e$ generated by
${\bar E}^\ell_\a$, ${\bar F}^\ell_\a$, $L^\ell_\mu$,
$\a\in\D_+$, $\mu\in P$.

The following is known (see  [DP] and references therein):
\begin{prop} \label{zroot}
\begin{itemize}
\item The subalgebra $Z^M_0$ is a Hopf
subalgebra in $\U^M_\e({\g}) $ .
\item It is also a Poisson subalgebra in $Z^M_\e$ and
together with the Hopf algbera structure is a Poisson-Hopf algebra.
\item $Z^M_\e$ is integrally closed
\item $Z^M_\e$ is  a free $Z^M_0$ module of the rank $\ell^r$.
\item $\U^M_\e({\g}) $ is finite-dimensional over $Z^M_0$
with \ $\dim_{Z^M_0}(\U^M_\e({\g}))=\ell^{\dim{\g}}$.
\item There is an isomorphism of Poisson Hopf algebras $Z^M_0\simeq C(G^*_M)$.
\end{itemize}
\end{prop}
Here $G_M^*$ is the finite covering of the
Poisson Lie group dual to the Poisson Lie group
$G$ (see section \ref{PL}). The isomorphism $Z^M_0\simeq C(G_M^*)$ is
given by the map $\phi$:
\[
\phi(L_\mu^\ell)=k_\mu, \ \phi(\bar{E}_i^\ell)=e_i, \ \phi(\bar{F}_i^\ell)=f_i.
\]

\begin{rem}
Geometrically, the algebra $\U^M_\e({\mathfrak g})$ can be regarded as a sheaf
of algebras over $G_M^*$ such that it is a bundle of algberas over each
symplectic leaf of $G_M^*$ with a flat connection over each simplectic leaf.
\end{rem}

\subsection{The completion $\overline{\U}^M_{\e}({\g})$ and
its center } Let $\ell$ be a positive integer.

\begin{prop} There exists a unique algebra structure
$\overline{\U}^M_{\!q}({\g})$ over ${\C}[q^{\pm 1}][[
q^{\ell}-1]]$ on the space formal power series
\[
\sum_{k_1\dots k_N\geq 0\atop {{{m_1\dots m_N\geq 0}\atop
{n_1\dots n_N\geq 0}}\atop {s\geq  0}}} C^S_{\{k\}\{m\}\{n\}}
{\bar E}^{k_1}_{\b_1}\dots {\bar E}^{k_N}_{\b_N} \ \prod^N_{i=1}
\big(L^\ell_{\mu_i}-1\big)^{n_i}\, {\bar F}^{m_N}_{\b_N}\dots
{\bar F}^{m_1}_{\b_1}
\]
such that restricted to polynomials in ${\bar E}_\beta,{\bar
F}_\beta$, $L_{\mu_i}$, $q^\ell-1$,$q^{\pm 1}$ it coincides with
$\U^M_q({\g})$. Here $C^S_{\{K\}\{m\}\{n\}}\in {\C}[L^{\pm
1}_{\mu_i}, q^{\pm 1}][[q^{\ell}-1]]$, \ $N\!=\!|\D_+|$, \
$\b_1,\dots ,\b_N$ is a convex ordering on $\D_+$, $\mu_1,\dots
,\mu_r$ are generators of $M$ and ${\bar E}_{\b_i}$ and ${\bar
F}_{\b_i}$ are as in (\ref{order}) .
\end{prop}

Specializing $q$ to $\e$ as in the previous section we obtain the
\ completion $\overline{\U}^M_{\e}({\g})$ of $\U^M_{\e}({\g})$.

The center $\overline{Z}^M_{\e}=Z(\overline{\U}^M_{\e}({\g}))$ has a
natural Poisson algebra structure.
  The following proposition is the formal  version
of the proposition \ref{zroot}.

\begin{prop}\label{FQUA}

\begin{itemize}
\item The subalgebra $\overline{Z}^M_0\subset \overline{Z}^M_{\e}$
generated by formal power series in\\${\bar E}^{\ell}_{\a}, {\bar
F}^{\ell}_{\a},L^{\ell}_{\mu}\!-\!1$ is a Hopf-Poisson subalgebra
in $\overline{\U}^M_{\e}({\g})$.
\item $\overline{Z}^M_0$ is isomorphic to the formal group $F(G_M^*)$.
\item $\overline{Z}^M_{\e}$ is a free $Z^M_0$-module of rank $\ell^r$.
\item $\overline{\U}^M_{\e}({\g})$ is f.d.~over
$\overline{Z}^M_0$ with {\rm{dim}}${}_{\overline{Z}^M_0}(
\overline{\U}^M_{\e}({\g}))= \ell^{\rm{dim}\,{\g}}$
\end{itemize}\end{prop}

Let us introduce formal variables $z_{\mu}$ as
$L^\ell_\mu =\e_\mu\exp(z_\mu)$
where $\e_\mu$ are elements of the finite order which generate
the group of automorphisms of the covering map $G_M^*\to G^*$ in a
neighborhood of 1 ( the same $\e_\mu$ that were used in section \ref{epsilon}.
Then ${\bar L}_\mu = L_\mu\exp\left(- \
\frac{z_{\mu}}{\ell}\right)$ are elements of finite order.

\subsection{ The universal $R$-matrix for $\overline{\U}^P_{\e}({\g})$}

 For each positive root $\b$ let $z_{\b},{\bar
E}_{\b},{\bar F}_{\b}$ be corresponding elements of
$\overline{\U}^P_{\e}({\g})$. Let
$\overline{\U}^P_{\e}({\g})^{{\hat\ot} 2}$ be the completion of the
tensor product with respect to the gradation given by the degree
function $deg(z_\b)=deg(\bar{E}_\b)= deg(\bar{F}_\b)=1$. Define
the outer automorphism ${\mathcal R}^{(n)}_{\b}$ of
$\overline{\U}^P_{\e}({\g})^{{\hat\ot} 2}$ as
\[
{\R}^{(n)}_\b(x) = \exp \left(\frac{1}{\ell^2} Li_2\big(  {\bar
E}_{\b}^{\ell}\ot {\bar F}_{\b}^{\ell}\big)\right)(x)
\]
where \ $\exp(y)\circ x=\ds{\sum^{\infty}_{n=0} \frac{1}{n!}}\,
\{y\{y\dots \{y,x\}\dots\}\}$ and
\[
Li_2(x)=-\int_0^x \frac{log(1-s)}{s}ds =\sum_{n=1}^\infty \frac{x^n}{n^2}
\]

Define the outer automorphism \ ${\mathcal R}^{(n)} =
{\mathcal R}^{(n)}_{\b_N}\circ\dots\circ {\mathcal R}^{(n)}_{\b_1}$.

Define the element \ ${\tilde R}^{(n)}_\b \in
\overline{\U}^M_{\e}({\g})^{{\hat\ot} 2}$ \ as
\[
{\tilde R}^{(n)}_\b =
\prod^{\ell -1}_{m=0} \big( 1-\e^m {\bar E}_\b \ot {\bar
F}_\b\big)^{-\frac{m}{\ell}} \ .
\]

Define the element \
$R^{(n)}\in \overline{\U}^P_{\e}({\g})^{{\hat\ot} 2}$ as
\[
R^{(n)}=R^{(n)}_{\b_N}\dots R^{(n)}_{\b_1}
\]
where
\[
R^{(n)}_{\b_i} = {\mathcal R}^{(n)}_{\b_1}\circ\dots\circ
{\mathcal R}^{(n)}_{\b_i-1}\big({\tilde R}^{(n)}_{\b_i}\big) \ .
\]

Define the outer automorphism ${\mathcal R}^{(c)}$ of
$\overline{\U}^P_{\e}({\g})^{{\hat\ot} 2}$
as
\[
{\mathcal R}^{(c)}(x)=\exp(\frac{1}{2\ell^2}\sum_{i,j=1}^r (b^{-1})_{ij}
z_{\alpha_i}\otimes z_{\alpha_j})(x)
\]
where $b_{ij}=a_{ij}d_j$ is the symmetrized Cartan matrix.

Define the element $R^{(c)}\in \overline{\U}^P_{\e}({\g})^{{\hat\ot} 2}$
as
\[
R^{(c)}=\sum_{\lambda, \mu \in P/dlP} \e^{(\lambda, \mu)}
P_\lambda\otimes P_\mu
\]
Here $d$ is the degree of the covering map $G^*_P\to G^*$ and
 $P_\lambda$ are idepompotents in the subalgebra generated by
$\bar{L}_\mu$ such that
\[
\bar{L}_\mu P_\lambda=P_\lambda\bar{L}_\mu =\e^{(\lambda,\mu)}
P_\lambda
\]

\begin{thm} \label{m-YB}\begin{enumerate}
\item The outer automorphism $\mathcal R$ restricted to
$\bar{Z}_0\hat{\otimes} \bar{Z}_0$ is a Poisson automorphism.
\item The automorphism ${\mathcal R}={\mathcal R}^{(c)}\circ
{\mathcal R}^{(n)}$ restricted to $\bar{Z}_0\hat{\otimes} \bar{Z}_0$
satisfies the Yang-Baxter equation
\begin{equation}\label{class-ybe}
{\mathcal R}_{12}\circ {\mathcal R}_{13}\circ {\mathcal R}_{23}
={\mathcal R}_{23}\circ {\mathcal R}_{13}\circ {\mathcal R}_{12}
\end{equation}
Here ${\mathcal R}_{ij}$ act on $\bar{Z}^{\hat{\otimes}3}$,
${\mathcal R}_{12}={\mathcal R}\otimes id$, ${\mathcal R}_{23}=
id\otimes {\mathcal R}$, ${\mathcal R}_{13}=id\otimes \sigma_{23}\circ
{\mathcal R}_{12}\circ id\otimes \sigma_{23}$ and $\sigma_{23}
(x\otimes y\otimes z)=x\otimes z\otimes y$.
\item The element $R=R^{(c)}R^{(n)}$ satisfies the twisted Yang-Baxter
relation
\begin{equation}\label{tw-YBE}
({\mathcal R}_{12}^{-1}\circ {\mathcal R}_{13})(R_{23})\cdot
 {\mathcal R}_{12}^{-1}(R_{13})\cdot R_{12} =
({\mathcal R}_{23}^{-1}\circ {\mathcal R}_{13}) (R_{12})
{\mathcal R}_{23}^{-1}(R_{13})R_{23} \ ,
\end{equation}
\item $\D'(a)={\mathcal R}(R\D(a)R^{-1})$ for all $a\in \bar{\U}_\e
({\mathfrak g})$.
\end{enumerate}
\end{thm}
This theorem follows from the asymptotical behavior of the universal
$R$-matrix for $sl_2$ \cite{R-2}, from the multiplicative formula
for the universal $R$-matrix and from teh Campbell-Hausdorf formula.
For arbitrary simple Lie algebra it was proven in \cite{Ga}.
The Poisson automorphism $\mathcal R$ was studied in \cite{WX} and
\cite{R-1}.

\section{Representations of $\U^P_\e({\g})$ and holonomy
$R$-matrices}

\subsection{Representations of $\U^M_\e({\g})$}

We will denote
by $(\pi_x^V,V)$ a representation $\pi_x^V: \U^M_\e({\g})\to End(V)$
of $\U^M_\e({\g})$ in the vector space $V$
with $Z^M_0$-central character $x\in G_M^*$ . Here we used
a natural identification of Poisson-Hopf algebras
$Z^M_0\simeq C(G^*_M)$ (see proposition \ref{FQUA}).

The group $G$ acts on $G_M^*$ locally by dressing transformations.
This  $G$-action on $G^*_M\simeq Spec(Z_0)$ lifts to the action
of $G$ on $\U^M_\e({\g})$ by outer
auto\mor s. We will write $g:a\mapsto \tau_g(a)$ for this action.

We will say that the representation $(\pi_x^V,V)$ is
{\it $G$-equivalent} to the representation $(\pi_y^W,W)$ if
\begin{itemize}
\item $x,y\in G_M^*$ belong
to the same $G$-orbit, i.e. if there exists $g\in G$ with $y=g(x)$
\item if there exists a linear map $\varphi_{V,W}(x,g):V\to W$ such that
\[
\pi^W_{g(x)}(a)=\varphi_{V,W}(x,g)\pi^V_x(\tau_g(a))\varphi_{V,W}(x,g)^{-1}.
\]
\end{itemize}

The representation of $\U^M_\e({\mathfrak g})$ dual to $(\pi_x^V, V)$
is the representation in the dual vector
space $V^*$ with the algebra acting as $a\mapsto \pi_x^V(S(a))^*$.
Here $S$ is the
antipode and $f^*$ is the linear map dual to $f:V\to V$.
It is clear that if $(\pi_x^V, V)$ has $Z^M_0$-character
$x\in G_M^*$ then the dual representation will have the
$Z^M_0$-central character $i(x)$ where $i$ is the operation of
taking the inverse in the group $G^*$.

\subsection{Irreducible representations}

Because the \alg \ $\U^M_\e({\g})$ is finite-dimensional over
its center, there exists a non-empty Zariski open
subset $S^M_\e\subset\Spec\!(Z^M_\e)$ such that
$\U^M_\e({\g})/\<(c-\chi(c))\U^M_\e({\g})\mid c\in Z^M_\e\>$ is
isomorphic to a matrix algebra for any $\chi \in S^M_\e$
( for more details see  \cite{DP}. We will call such
elements generic. Thus, for each generic  $\chi\in \Spec\!(Z^M_\e)$
we have unique iso\mor \
class of irreducible representations.

Denote by $S^M_0\subset Spec(Z^M_0)$ the image of
$S^M_\e$ with respect to the projection $Spec(Z^M_\e)\to Spec(Z^M_0)$.
The variety $S^M_\e$ is a finite cover of $S^M_0$. Points of
$Spec(Z^M_\e)$ are "common level surfaces" of all central elements
of $\U^M_\e({\g})$. Points of $Spec(Z^M_0)$ are "common level surfaces" of
elements of central subalgebra generated by $L_\mu^\ell, \bar{F}_i^\ell,
\bar{E}_i^\ell$ and by their Poisson brackets. The number of branches
of the projection $S^M_\e\to S^M_0$ over generic point is $\ell^r$.
The number of branches of $S^M_0\to S^Q_0$ is $\ell^{|M/Q|}$.
So, over generic point of $G^*$ we expect $\ell^{r+|M/Q|}$ irreducible
representations of $\U^M_\e({\g})$. All these irreducible representations
have dimension $\ell^{|\D_+|}$ \cite{DKP}.

Central elements of $\U^M_\e({\mathfrak g})$ which also
belong to the Poisson center of $Z^M_\e$
are constant on dressing
orbits. We will
call this central subalgebra the Casimir subalgebra.

Let ${\mathcal O}\subset G^*$ be a dressing orbit and $U\in G^*$ be a
neighborhood of $1$ on which the local action of $G$ integrates
to an action. Let $\{(\pi_x^V,V)|x\in U\subset {\mathcal O}\}$ be a family of
irreducible representations of $\U^M_\e({\g})$. Assume that these
representations have the same central character with  respect to
the Casimir subalgebra. Because Poisson Casimirs are constant on $G$-orbits and
the specter of primitive ideals is a finite cover over the specter of
primitive ideals of the Poisson center of $Z^M_0$, this assumption will hold for
sufficiently small $U$. Let $g\in G$ and $x\in
{\mathcal O}\cap U\subset G^*$ be such that $x$ and $g(x)$ are generic.
Representations $\pi_x^V\circ \tau_g$ and $\pi_{g(x)}$ have the same
central characters and therefore isomorphic. Thus, we have a family $\{T(g|x)\}$ of
linear automorphisms of $V$ such that
\begin{equation}\label{connect}
\pi_x^V(\tau_g(a))=T(g|x)\pi_{g(x)}^V(a)T(g|x)^{-1}
\end{equation}

Considering formal neighborhood of $1$ in $G$ and $G^*$ and
representations of $\U^M_\e(\g)$ over such neighborhood we
get what is called formal representations. Such representations are
homomorphisms from $\bar{\U}^M_\e(\g)$ to the algebra $End(V)[[x]]$
where $V$ is the space of representations and $x$ are formal coordinates
in a neighborhood of $1$.

\subsection{Holonomy $R$-matrices}

From now on we will use the map $I$
to identify neighborhoods of identities in $G^*$ and $G$.
After this the dressing action of $G$ on $G^*$ is identified with the
adjoint action of $G$ on itself.

Let $(\pi_x^VV)$ and $(\pi_y^W,W)$ be two generic formal
representations of $\U^P_\e({\mathfrak g})$.
\begin{prop}\label{eval-R} Let $a, b\in \U^P_\e({\mathfrak g})$ and $\mathcal R$ be
the outer automorphism of $\U^P_\e({\mathfrak g})^{\otimes 2}$ defined in
the theorem \ref{m-YB}. Then
\begin{equation}
(\pi_x^V\otimes \pi_y^W)({\mathcal R}(a\otimes b))=
\pi_x^V(\tau_{x_L(x,y)_+}(a))\otimes \pi_y^W(\tau^{-1}_{x_-}(b))
\end{equation}
where $x_L(x,y)=x_-yx_-^{-1}$
\end{prop}

Proof. This proposition follows from the definition of $\mathcal R$
and from results of \cite{WX}.

For generic formal $x$ and $y$ define the element $R^{V, W}(x,y)\in
End(V\otimes W)[[x,y]]$ as
\begin{equation}\label{hol-R}
R^{V,W}(x,y)=(T(x_L(x,y)_+|x)^{-1}\otimes T(x_-|y)^{-1})(\pi_x^V\otimes
\pi_y^W)(R)
\end{equation}
Here $x$ and $y$ are formal coordinates on a formal neighborhood of $1$ in
$G$.

\begin{thm} Linear maps (\ref{hol-R}) satisfy the holonomy Yang-Baxter
equation:
\begin{eqnarray}\label{hYBE}
R^{V,W}_{12}(x_R(x,x_L(y,z)),x_R(y,z)) R^{V,U}_{13}(x,x_L(y,z))
R^{W,U}_{23}(y,z)=
\\ R^{W,U}_{23}(x_L(x,y),x_L(x_R(x,y),z)) R^{V,U}_{13}(x_R(x,y),z) R^{V,W}_{12}(x,y)
\end{eqnarray}

\end{thm}

Proof. We can choose linear isomorphisms $T(g|x)$ such that
$T(g_1g_2|x)$ is proportional to $T(g_1|x)T(g_2|x)$ and
\begin{equation}\label{T-invert}
T(g^{-1}|x)T(g|gxg^{-1})=1
\end{equation}

From Proposition \ref{eval-R} we can evaluate the both sides of the
equation (\ref{tw-YBE}) in the tensor product of three $G^*_P$ evaluation representations.
For the left side we have:
\begin{eqnarray}
(\pi_x^V\otimes\pi_y^W\otimes\pi_z^U)(({\mathcal R}_{12}^{-1}\circ
{\mathcal R}_{13})(R_{23})\cdot
 {\mathcal R}_{12}^{-1}(R_{13})\cdot R_{12} )=\\
 (\pi_x^V\otimes\pi_y^W\otimes\pi_z^U)(\tau_{x_L(x,x_L(y,z)_-)_+}\otimes\tau_{x_L(y,z)_+})(R)_{12}
 (id\otimes\tau_{y_-})(R)_{13}R_{23}
 \end{eqnarray}
Similarly one can evaluate the right side:
\begin{eqnarray}
(\pi_x^V\otimes\pi_y^W\otimes\pi_z^U)(({\mathcal R}_{23}^{-1}\circ {\mathcal R}_{13}) (R_{12})
{\mathcal R}_{23}^{-1}(R_{13})R_{23} )=\\
(\tau_{x_-}^{-1}\otimes\tau_{x_R(x,y)_-}^{-1})(R)_{23}(\tau_{x_L(x,y)_+}^{-1} \otimes
id)(R)_{13}R_{12}
\end{eqnarray}
where $x_R(x,y)=x_L(x,y)_+^{-1}xx_L(x,y)_+$.

The holonomy Yang-Baxter equation for linear maps (\ref{hol-R}) follows
from the identities (\ref{connect}) and the identities for $T(x|y)$.

Thus we constructed solutions to formal holonomy Yang-Baxter equation.
Let $S$ and $S'$ be two generic symplectic leaves in $G^*_P$.
Linear operators (\ref{hol-R}) admit analytical continuation to sections of
vector bundles over $S\times S'$.

\section{Invariants of tangles}
\subsection{d-matrix}

In the construction of invariants of tangles with flat connection
given in \cite{KR} an important role played linear operators
$d^V(x)\in End(V)$ defined in terms of holonomy $R$-matrices as
\[
d^X(a)=(tr\otimes id)(P((R^{t_1}(a,i(a)^{-1})^{-1})^{t_1})
\]
Since we constructed holonomy $R$-matrices for irreducible
$\U^P_\e(\g)$-modules, we have such $d$ operators for each
generic irreducible representation $(\pi_x^V,V)$.

\begin{lemma}Let $(\pi^V_x, V)$ be an irreducible representation of
$\U^P_\e(\g)$ with generic $x\in G^*$, then
\[
d^V(x)=c_V(x)\pi^V_x(L_\rho)
\]
where $c_V(x)$ is a non-zero complex number
and $\rho=1/2\sum_{\alpha\in \D_+}\alpha$.
\end{lemma}

Proof. For each pair of representations $(\pi^V_x, V)$ and $(\pi^W_x, W)$
of $\U^P_\e(\g)$ we have:
\begin{equation}\label{cr-1}
R_{12}^{X,Y}(x,y)=d_2(a)^{-1}(((R_{12}^{X,Y}(x,y)
^{-1})^{t_2})^{-1})^{t_2})d_2(y)
\end{equation}
\begin{equation}\label{cr-2}
R_{12}^{X,Y}(a,b)=d_1(y)^{-1}(((R_{12}^{X,Y}(a,b)
^{t_1})^{-1})^{t_1})^{-1})d_1(a)
\end{equation}
These equations imply that for generic $x\in G^*$ and a representation
$(\pi^V, V)$ we have:
\[
\pi^V_x(S^2(a))=d^V(x)\pi^V_x(a)d^V(x)^{-1}
\]
On the other hand from the definition of the antipode it is
easy to see that
\[
S^2(a)=L_\rho aL_\rho^{-1}
\]
The lemma now follows from the Schur's lemma.

This lemma implies that
the invariant of a knot defined by
the functor $F$ constructed in \cite{KR} is identically zero.
The situation is similar to invariants studied in \cite{Ro} (see also
the references therein).

\subsection{Invariants knots of string knots}

Recall that a string knot is a tangle with one connected component
and with two boundary components. If $D_t$ is a diagram of a string knot
a generic $G$-coloring of the whole diagram is determined by the
corresponding  $G$-coloring of one of its boundary component.

\begin{prop}
Let $t$ is a string knot and $F_V(t,x)$ is the value of the functor $F$
on it. Here we assume that the lower and upper boundaries of $t$ are $G$-colored by
$x$ and decorated by $\U^M_\e({\g})$-module $(\pi_x^V, V)$. Then the
element $F_V(t,x)\in End(V)$ is invariant with respect to gauge
transformations, i.e. $F_V(t^g,g(x))=F_V(t,x)$ where $t^g$ is the
result of the gauge action of $g\in G$ (see \cite{KR}) on the $G$-colored tangle $t$
and $g(x)$ is the result of the dressing action of $g\in G$ on $x$.
\end{prop}
Proof. The gauge invariance of the functor $F$ (see \cite{KR}) implies
that for any other representation $(\pi^W_x, W)$ of $\U^M_\e({\g})$
we have:
\[
(1\otimes F_V(t^{y_+},y_+(x)))R^{V,W}_{21}(x_R(y,x),x_L(y,x))^{-1}=
R^{V,W}_{21}(x_R(y,x),x_L(y,x))^{-1}(1\otimes F_V(t,x))
\]
\[
(1\otimes F_V(t^{y_-},y_-(x)))R^{W,V}_{12}(y,y_-(x))^{-1}=
R^{W,V}_{21}(y,y_-(x))^{-1}(1\otimes F_V(t,x))
\]
These equations with $y=1$ imply that elements $F_V(t,x)$ are central.
Same equations for generic $y\in G$ imply that
$F_V(t^y,y(x))=F_V(t,x)$.

For generic $q$ the functor $F$ (see \cite{RT}) applied to a string
knot defines an element of a completion of the center of $\U_q(\g)$.
This element can be evaluated in a finite-dimensional representation
and up to a scalar factor (which is equal to the quantum dimension
of the representation) coincides with the corresponding invariant
of the knot obtained by closing the string knot. This means that
the central element itself is not only an invariant of a string knot but is
also an invariant of the knot obtained by closing the string knot.

One can argue that the same happens in our case. The "limit" of this
central element when $q\to \e$ according to the asymptotical behavior
of the universal $R$-matrix \cite{R-2} has an essential singularity and
and a regular part. The essential singularity will give the invariant
related to the Poisson $R$-matrix ( see \cite{WX} and \cite{R-1}).
The regular part will give the invariant discussed here.
As it was explained above for generic $q$ the
central element which is the invariant
of a string knot is also an invariant of the knot obtained by
the closing the string knot. Therefore, we should expect
the same for roots of 1.
We will return to the detailed discussion of this question
in a separate publication.

\section{Conclusion}

We constructed invariants of tangles with flat connections
over the complement using representation theory of
quantized universal enveloping algebras at roots of 1.

We conjecture that these invariants for $SL_2$ coincide with the
invariants constructed in \cite{BB} in case when the $3$-manifold is a
complement to a tangle. When $G=SL_2$ and the flat connection is trivial
they coincide with the invariant constructed in \cite{Ka}.

Since the invariant for $G=SL_2$ and trivial flat connection in the
complement gives the hyperbolic volume of the complement when $l\to
\infty$. It would be very interesting to describe the asymptotic of our new
invariants in terms of corresponding  geometrical invariants.

Let $\phi$ be a flat connection in the complement to a tangle.
We expect that in the limit $\phi\to 1$, where $1$ is the trivial flat
connection, our new invariant becomes the invariant constructed in
\cite{RT} for roots of 1 and reducible but indecomposable representations
of dimension $\ell^{\D_+}$. For $SL_2$ this gives the relation between
the invariant constructed in \cite{Ka} and the Jones polynomial at roots
of 1, which was observed in \cite{MM}.

What has been done in our two papers is a first step in the larger
program. Here we will outline of some further steps.

First question is whether there is a topological quantum field theory
which can give a geometric description of these invariants.
In case of Jones polynomials such phenomenological quantum field theory
(Chern-Simons theory for compact simple Lie groups) was proposed by
Witten and allowed to describe the invariants in geometrical terms.
One can guess that appropriate version of Chen-Simons theory
for complex simple Lie groups will describe the large $l$ asymptotic
of our invariants.

There is another description of our invariant in terms of triangulated
manifolds which is based on "6j-symbols" for the category of modules
over $\U_\e(\g))$. It generalizes the construction from \cite{TV}\cite{Ka}
and \cite{BB}. This construction also gives invariants
of more general 3-manifolds with flat connections. We will do it in
a separate publication.

\end{document}